\documentclass[journal, onecolumn, dvipdfm]{IEEEtran}
\usepackage[justification=centering]{caption}
\usepackage{cite}
\usepackage[pdftex]{graphicx}
\usepackage{amsmath}
\usepackage{amssymb}
\usepackage{url}
\usepackage{setspace}
\renewcommand{\baselinestretch}{1.5}
\usepackage{authblk}
\begin{document}

\title{Extended SPR for Evolving Networks with Nodes Preferential Deletion}

\author{
\IEEEauthorblockN{Yue~Xiao\IEEEauthorrefmark{2}, Xiaojun~Zhang\thanks{* Xiaojun~Zhang is the corresponding author.}\IEEEauthorrefmark{3}}\\
\IEEEauthorblockA{School of Mathematical Science, University of Electronic Science and Technology of China\\ Chengdu, Sichuan 611731, China
}
\\
\IEEEauthorblockA{\IEEEauthorrefmark{2}yuexiao@std.uestc.edu.cn \\
\IEEEauthorrefmark{3}sczhxj@uestc.edu.cn
}
}

\maketitle

\begin{abstract}
Evolving networks are more widely existed in real world than static networks, and studying their statistical characteristics is vital to recognize and explore them further. But for the networks with nodes preferential deletion, there are few researches due to the lack of effective methods. In this article, we propose an extended SPR (ESPR) for these preferential removal networks when discuss the essential statistics, especially the steady-state degree distribution. Comparing with continuum formalism that is often employed, this theory-supported method retains the actual topological structure and corresponding statistics of networks during evolving process. With two theorems, we demonstrate the effectiveness of ESPR in handling evolving networks with nodes non-uniform removal; moreover, it also be proved that the SPR is special case of ESPR. In other words, ESPR is an operative framework when deal with the degree distibution, and it even have potential to solve other statistics of evolving networks. 
\end{abstract}


\IEEEpeerreviewmaketitle

\textbf{\section{Introduction}}
\indent\IEEEPARstart{C}{omplex} networks are a special class of intuitive structures that involves a large amount of connected nodes with edges, and it can be extracted from many fields, like social networks\cite{newman2003random}, biological networks\cite{barabasi2004network}, and power grids\cite{wang2009cascade}. To recognize topological structure and explore dynamic behaviors of complex networks, researchers have focused on statistical characteristics of them, such as degree distribution, degree-degree correlation etc. with the development of related theory for years\cite{watts1998collective, barabasi1999emergence, albert2002statistical, newman2003structure}. In fact, different from the static networks data collected deliberately, there are more evolving networks in reality, in which nodes and edges attached on them would be added or removed as time passes. For example, terminals with less neighbors would be easier to be removed from the internet; and people involved in social network tend to contact new friends  continuously\cite{dorogovtsev2002evolution, slater2006evolving,ferretti2011preferential,he2007time}. It is hard to discuss these complicated statistics of dynamic networks accompanying different evolving mechanisms directly.\\
\indent Among these essential characteristics, the steady-state degree distribution has been discussed since the introduction of BA model and further study of network scaling \cite{li2003local}. For this discrete statistics, researchers have tried to solve it with continuum equaitons when process the pure growing network models. Specifically, Barab{\'a}si and Albert proposed the mean-field equation to solve the steady degree distribution of the BA model. They assumed that the degree $k$ is continuous, the probability of new node connects to existing nodes in the network can be interpreted as a change rate of $k$\cite{barabasi1999emergence}. Similarly, Dorogovtsev proposed the master-equation approach for a model with preferential linking and all nodes have initial attractiveness. Basing on the distribution $p(q,s,t)$ of the connectivity $q$ of the node $s$ at $t$, they even discussed the relationship between the exponent of power-law degree distribution and different values of parameters\cite{dorogovtsev2000structure}.  And Krapivsky et al. also applied the rate-equation approach to study the time and age dependent connectivity distribution of a growing random network\cite{krapivsky2000connectivity}. Besides, different from former continuum approaches, Shi et al. established the relation between growing networks and Markov chains firstly and calculated the degree distribution with a computational approach they established\cite{shi2005markov}. These methods and their variants are also applied to solve the growing networks which remove nodes, especially in the case of nodes random deletion\cite{moore2006exact,ben2007addition,Salda2007Continuum,garcia2008degree,2008Degree, ikeda2019graph}. Typically, specific continuous approach was proposed to solve models with nodes random removal and discuss the scaling properties of its scale-free degree distribution\cite{dorogovtsev2001scaling}. Other scholars also tried to discuss the degree distribution of this kind of evolving networks with different strategies from a queueing system perspective\cite{0Heritable}, or using the relative entropy\cite{2020Analysis}.\\
\indent According to the researches introduced above, we can conclude that continuum formalism is indeed convenient in processing this discrete problems, and we can obtain the approximate results for the steady-state degree distribution. But these continuous approaches only consider the average change of degree for single node without the actual topological structure of networks, and the fact is that degree distribution is a vital statistic reflects the global feature of networks. Besides, there are also various possible networks appearing with different probabilities in the evolving process, and each of them own steady-state degree distribution. Moreover, comparing with those networks introduced above, another removal mechanism of evolving networks, nodes preferential deletion has fewer studies, which can be attributed to the lack of appropriate method to some extent. Although some researches focused on the stability, convergence and even uncorrelatedness of these networks with nodes removed non-uniformly\cite{2004Scale,2005Stability,2007A,2011L,2014Tracking, 2020Analysis}, how to calculate the steady-state degree distribution for corresponding evolving networks is quite difficult. Only Kong et al. discussed this problem with a mean-field rate equation approach\cite{2008Preferential}. Deserve to be mentioned that  Zhang et al. proposed a stochastic process rule (SPR) considering the special structure of networks, it focuses on the transfer approaches of a node with two-dimensional state $(n,k)$ in a inifinte system $\Omega$, where $n$ is the scale of network the node located and $k$ is the degree of the node. So we can describe the evolving process with a transition equation $P(t+1)=P(t)P$ where the $P(t)$ and $P(t+1)$ are the probability distributions of networks at $t$ and $t+1$ respectively, $P$ is the one-step transition probability matrix. SPR is an efficient method that remains the same topologies and statistics as generated networks the next moment. As a theory-supported method, SPR has been used to solve the BA model with a result of power-law degree distribution\cite{2012SPR}, and a random birth and death network(RBDN) has also been studied successfully\cite{zhang2016random}. Therefore, it seems to be possible to discuss evolving networks with nodes preferential removal with principle of SPR. Suppose a network decays at $t$, each node in this network would be removed with a probability proportional to its degree $k$ or these probabilities depend on the degree distribution of network, not only the number of nodes in the network. The one-step transition equation can be used to illustrat this process, though the form of transition probability would be more complicated than that of nodes random deletion.\\
\indent In this paper, we promote the original SPR for evolving networks with nodes preferential deletion in handing the steady-state degree distribution. It is organized as follows: in Section \ref{sec.2}, we introduce the essential definitions firstly; in Section \ref{sec.3}, we interpret the ESPR for nodes preferential deletion basing on the principle of SPR; in section \ref{sec.4}, we proposed two theorems to confirm the effectiveness of ESPR; conclusions and prospects of future work based on ESPR are illustrated in Section \ref{sec.5}.\\

\section{Preliminaries}\label{sec.2}
\indent\textbf{Definition1}.$(Degree\quad Distribution)$ 
For graph $G=<V,E>$, where $V$ is the set of nodes and $E$ is the set of edges, let $K_G$ be a random variable, then
\begin{equation}
	P\left\{K_G=k\right\}=\frac{N_k}{N}
\end{equation}
is the degree distribution of $G$, where $k=0,1,2,...,|V|$ and $N_k$ is the number of nodes have $k$ edges and $N=|V|$.

Further, for a graph process $\left\{G(t),t\geq 0\right\}$, a series of graphs would be generated with different probabilities between $t$ and $t+1$ instead of an average one. The degree distribution for a graph mentioned above is limited to illustrate statistics of this process. Thus, we give the definition of average degree distribution for $\left\{G(t)\right\}$ is shown as below.\\

\indent\textbf{Definition2}.$(Average\quad Degree\quad Distribution)$
For graph process $\{G(t),t\geq 0\}$with $P_{G_i}(t)=P\left\{G(t)=G_i(t)\right\}$, where $G_i(t)$ represents these networks which own $i(i=1,2,...,N(t))$ nodes, $N(t)$ is the largest network in scale at $t$, let $K(t)$ be a random variable, then the average degree distribution $P_k(t)$ of $\left\{G(t)\right\}$ at $t$ is defined as 
\begin{equation}
	\begin{aligned}
		P_k(t)&=P\left\{K(t)=k\right\}\\
		&=\sum_{i=1}^{N_{t}}P\left\{G(t)=G_i(t)\right\}P\left\{K_{G_i(t)}=k\right\}\\
		&=\sum_{i=1}^{N_{t}}P_i(t)P\left\{K_{G_i(t)}=k\right\}.
	\end{aligned}
 \end{equation}

According to the definition of degree distribution in\cite{dorogovtsev2000structure}, here, we define the steady degree distribution for $\left\{G(t), t\geq 0\right\}$ as follows.\\

\indent\textbf{Definition3}.$(Steady-state\quad Degree\quad Distribution)$
For graph process $\left\{G(t), t\geq 0\right\}$, and $K(t)$ is average degree of it at $t$. Let $K$ be a random variable, then the steady degree distribution $P_k$ of $\left\{G(t)\right\}$ is defined as
\begin{equation}
\begin{aligned}
P_k&=P\left\{K=k\right\}\\
&=lim_{t \rightarrow + \infty}P\left\{K(t)=k\right\}. 	
\end{aligned}
\end{equation} \\

\section{ESPR for Nodes Preferential Deletion}\label{sec.3} 
\subsection{SPR}
For an evolving network, adding and removing nodes would influence its topological structure, and this process cannot be described directly with continuum approaches. However, we can illustrate the evolving process of networks from the perspective of nodes with SPR. Specifically, an infinite sample space $\Omega$ for networks generated during evolving process $\left\{G(t), t\geq0\right\}$ is introduced, $G(0)$ is the initial network which is usually a compelete graph by default. $G(t)=\cup_{i}^{N(t)}G_i(t)$ is a division of $\Omega$ at $t$, where $G_i(t)$ includes all networks that own $i$ nodes though they may be different in topology at $t$, and $N(t)$ is the upper bound of $i$ at $t$. \\
\indent Basing on the particular structure of network, arbitrary node in a chosen network would possess a two-demensional state $(n,k)$, where $k$ is the degree of this node and $n$ is the scale of network on which it is located. Similarly, two nodes in $\Omega$ may have the same state here, though they may be completely different from the perspective of the time they are added into the network. Then the change of network either growing or decaying from $t$ to $t+1$ could be described by the transition of nodes state with previous setting, e.g. a node with state $(n,k)$ at $t$ would be $(n+1,k)$ or $(n+1,k+1)$ at $t+1$ when network grows with probability $p$, and may also be $(n-1,k)$ or $(n-1,k-1)$ at $t+1$ when network decays with probability $q=1-p$. Two schematic diagrams for evolving networks with adding and deleting nodes randomly are given in Appendix respectively, showing the specific transition process of nodes in SPR. In diagrams, we give these nodes different numeric notations from $1$ to $N(t)$ to distinguish each other in a network, and some of them with different notations may still have the same state. Additionally, the $q_v$ is the removing probability for node denoted as $v(v=1,2,...,N(t))$, $P_{(n,k)(n^{'},k^{'})}(t)$ is the transition probability of node $v$ from state $(n,k)$ to $(n^{'},k^{'})$ at $t$, $k^{'}$ is probable degree of $v$ at $t+1$. \\
\indent According to the principle and related definitions above, the framework for solving steady-state degree distribution with SPR is given as follows:\\

		
\indent$(\romannumeral1)$ For a node $v$, determine the state transition probability $P_{(n,k)(n^{'},k^{'})}(t)$;\\
\indent$(\romannumeral2)$ Give the complete state transition equation of $v$ from $t$ to $t+1$;\\
\indent$(\romannumeral3)$ Determine the transition equations of degree distribution $P_k(t)$ for different degree $k(k=0,1,2,...)$;\\
\indent$(\romannumeral4)$ Solve the corresponding equations of steady degree distribution $P_k$ for different $k$.\\

Among these progressive steps, the first step is the most significant one absolutely, which makes the SPR distinguished from continuum formalism. And it is worth mentioning that two degree distribution in stationary from two model still be the same, once they had the same state transition equation for considered node\cite{2012SPR}. So proofs that in section \ref{sec.4} only need to compare the results of the first two steps.\\

\subsection{Extended SPR (ESPR)}

According to the principle shown in Fig.4, those isolated nodes deleted randomly from network at $t$ would be used to construct networks that appear at $t+1$ and corresponding state-transition probabilities can be given intuitively. But for nodes preferential deletion, original approaches mentioned above would be inadequate. Because in Fig.\ref{fig.4}, isolated nodes with different denotions are uniform from $t$ to $t+1$, which makes constructing new networks with isolated nodes at $t+1$ possible. However, if nodes in network are removed with probabilties propotional to their degrees respectively at $t$, the ratios of isolated nodes in different denotions depends on $t$ actually. Basing on the advantages of SPR and approaches for node random deletion, we try to give some improvements to handel nodes preferential deletion.\\
\indent Converting those nodes which are abundant in quantity to absent ones respectively and proceeding the reconstructure step may be operable. But it is quite difficult to give out state-transition probabilities. Here, we find a more appropriate operation which would been proved to unify both nodes random and preferential deletion in the form of state transition probability. Specific approaches corresponding to the first step in original SPR from $t$ to $t+1$ are given as follows:\\


\indent\textbf{Step1}: Select $\alpha(t)$ networks possessing the same topological structure from $G_n(t)(n=1,...i...,N(t))$, then choose one of them randomly and remove a node with a probability proportional to its degree. For two temporary sets $A(t)$ and $B(t)$, the newly generated network with $n-1$ nodes after node deletion and the removed node are placed into them respectively. Repeat the same  operations on the rest $\alpha(t)-1$ networks, and we will get $B(t)=\cup_vB_v(t)$ finally, where $B_v(t)$ represents the set of nodes with denotion of $v$ at $t$;\\
\indent\textbf{Step2}: For each kind of node denotions in $B(t)$, convert corresponding nodes to $N(t)$ types respectively with the ratios of that in $A(t)$, and put these newly generated nodes into the set $C(t)$;\\
\indent\textbf{Step3}: Referring to the denotions of nodes in a network selected from $A(t)$ randomly, choose $n-1$ isolated nodes from $C(t)$ with the same combination of types. And connect these isolated nodes with edges to achieve the same topology as the selected network, then put the newly generated network into a set $G_{n-1}(t+1)$. Repeat the same operations on the rest elements in $A(t)$ and $C(t)$, and the final $G_{n-1}(t+1)$ involves $n-1$ new networks at $t+1$ through evolving process.\\

\begin{figure}[!t]
	\centering
	\includegraphics[width=5.5in]{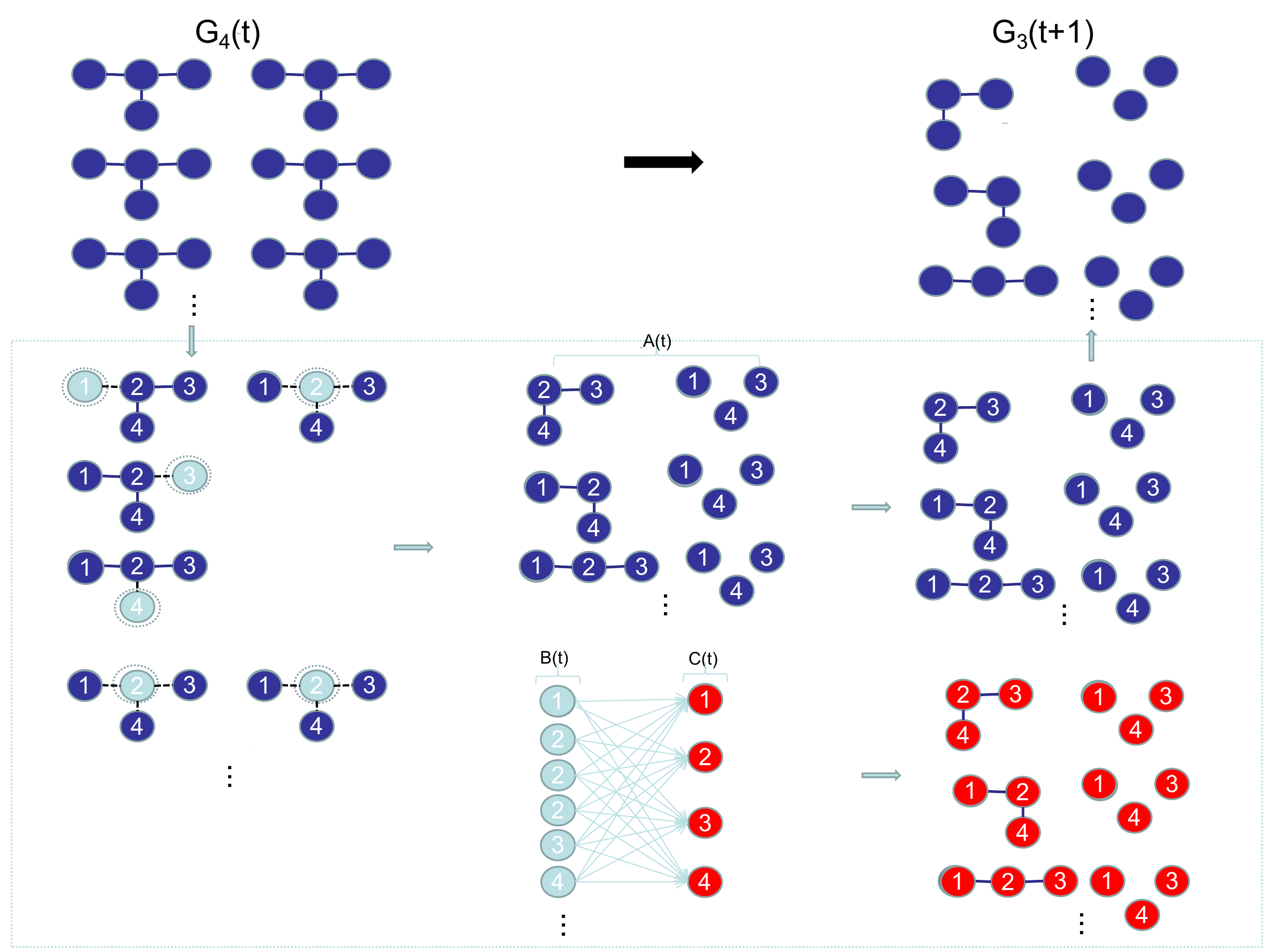}
	\caption{nodes preferential deletion.}
	\label{fig.1}
\end{figure}

Where the $\alpha(t)$ is a variable depends on the scale of network $n$. The schematic diagram of nodes preferential deletion is given as follow. Comparing the approaches shown in Fig.\ref{fig.1} with that in Fig.\ref{fig.4}, we can find that they are quite similar except the procedure to reconfigure isolated nodes, demonstrating that the improvement keeps the framework of SPR. \\ 

\subsection{State Transition Probability}
For a network having $n$ nodes at time $t$ and a selected node $v(v=1,2,...,n)$ from it, let $(n,k)$ be the state of $v$, and $N_v$ represents the set of  neighborhoods of $v$. Besides, we assume that $k^{'}$ represents probable degree of $v$ at $t+1$ with ESPR, $P_{(n,k^{'})}(t)$ is the probability of nodes with degree $k^{'}$ in the network at $t$, and $q_v(t)$ and $q_w(t)$is the removal probability of $v$ and another node $w$ at $t$ respectively. \\
\indent With ESPR, $v$ may be involved in $A(t)$ or $C(t)$ during transition, so the one-step transition probabilities of $v$ from $t$ to $t+1$ are shown as follows. \\
$(1)$ Node $v$ with state $(n,k)$ is not isolated:\\
\begin{equation}
\begin{aligned}
		P_{(n,k)(n-1,k)}(t)&=\sum_{w\notin N_v}q_w(t) \quad\quad\quad
		P_{(n,k)(n-1,k-1)}(t)&=\sum_{w\in N_v}q_w(t)\\
\end{aligned}
\end{equation}
$(2)$ Node $v$ with state $(n,k)$ is isolated:\\
\begin{equation}
\begin{aligned}
		P_{(n,k)(n-1,k^{'})}^{*}(t)&=q_v(t)[\frac{(1-q_{(n,k^{'})}(t))P_{(n,k^{'})}(t)}{\sum_jP_{(n,j)}(t)(1-q_{(n,j)}(t))}\frac{P_{(n,k^{'})(n-1,k^{'})}(t)}{P_{(n,k^{'})(n-1,k^{'})}(t)+P_{(n,k^{'})(n-1,k^{'}-1)}(t)} \\
		&\quad +\frac{(1-q_{(n,k^{'}+1)}(t))P_{(n,k^{'}+1)}(t)}{\sum_j(1-q_{(n,j)}(t))P_{(n,j)}(t)}\frac{P_{(n,k^{'}+1)(n-1,k^{'})}(t)}{P_{(n,k^{'}+1)(n-1,k^{'}+1)}(t)+P_{(n,k^{'}+1)(n-1,k^{'})}(t)}] \\
		&\quad =\frac{nq_{v}(t)}{n-1}[P_{(n,k^{'})}(t)P_{(n,k^{'})(n-1,k^{'})}(t)+P_{(n,k^{'}+1)}(t)P_{(n,k^{'}+1)(n-1,k^{'})}(t)]
\end{aligned}
\end{equation}
Note that we distinguish the transition probability of $v$ with $P^{*}$ when it be isolated during evolving, and $q_v(t)=q_{(n,k)(t)}$. \\
Therefore, for any node in a network with the state of $(n,k)$ at $t+1$, it may be transformed from other states and the state transition equation is illustrated as follows. 
\begin{equation}
\begin{aligned}
		P_{(n,k)}(t+1)&=P_{(n-1,k)}(t)P_{(n-1,k)(n,k)}(t)+P_{(n-1,k-1)}(t)P_{(n-1,k-1)(n,k)}(t)+P_{(n+1,k)}(t)P_{(n+1,k)(n,k)}(t)\\
		&\quad +P_{(n+1,k+1)}(t)P_{(n+1,k+1)(n,k)}(t)+\sum_{k^{'}}{P_{(n+1,k^{'})}(t)P_{(n+1,k^{'})(n,k)}^{*}(t)}
\end{aligned}
\end{equation}\\

\section{Theorems} \label{sec.4}
As we know, evolving networks in real world are obtained directly with setted evolutionary rules (ER), i.e., nodes are added into the network or removed from it at each moment, and only a single network sample is available at the next moment, rather than the whole sample space supposed in ESPR. Therefore, it is necessary to compare the average degree distributions obtained with both methods. If they are the same, it is feasible to discuss the statistics of the evolving network with ESPR.  \\

\begin{figure}[!t]
	\centering
	\includegraphics[width=5in]{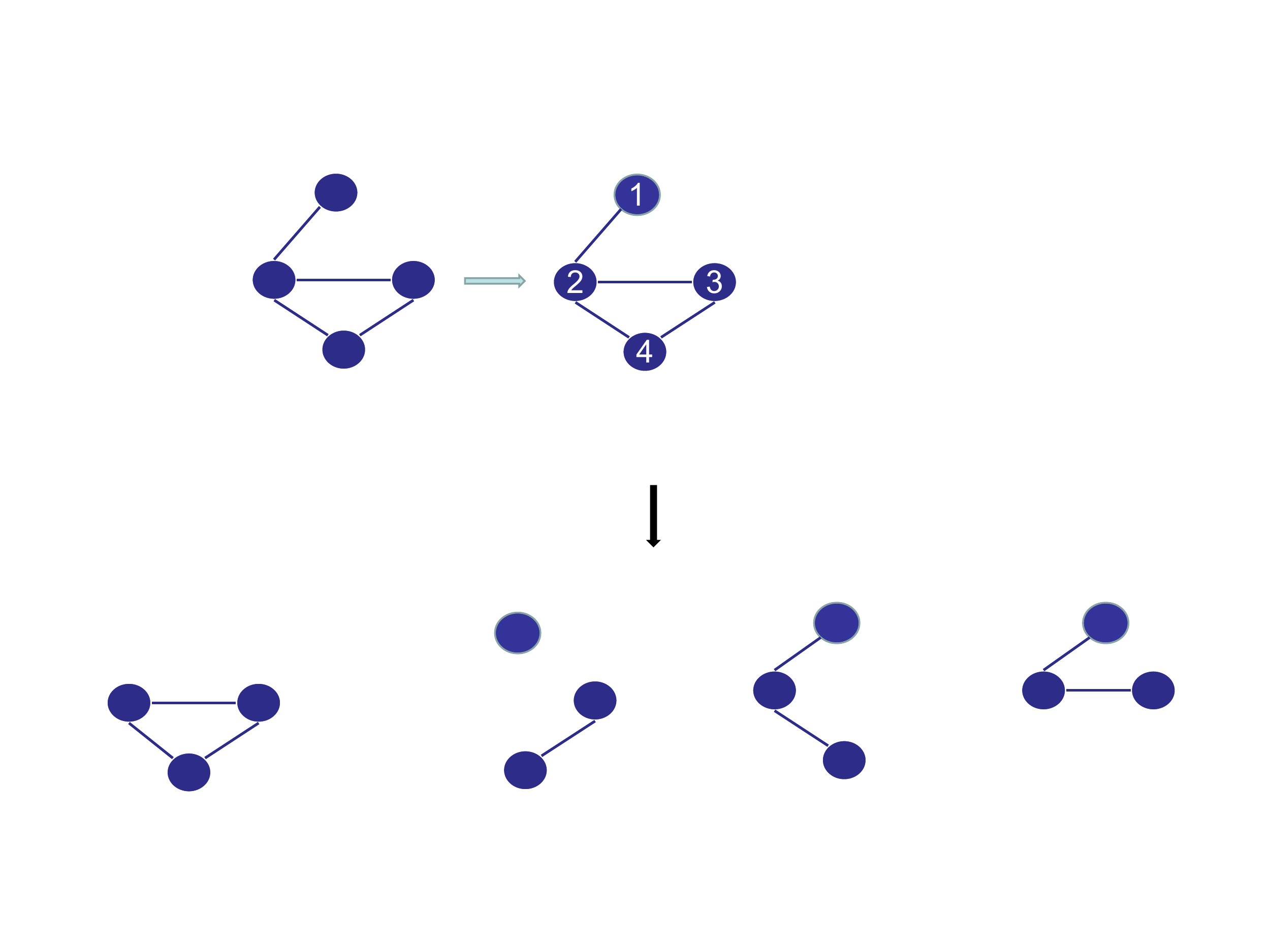}
	\caption{An schematic diagram of sample network.}
	\label{fig.2}
\end{figure}

\textbf{Theorem 1}. For an evolving model with nodes preferential deletion $\{G_{t},t\geq0\}$, suppose $EP_{k}$ represents the steady degree distribution of $\left\{G_{t}\right\}$ obtained from ESPR, and $P_{k}$ represents the steady degree distribution of it obtained from evolving rule directly. Then, the formula $P_{k}=EP_{k}$ hold. \\
\textbf{Proof}. The evolving process in reality is limited by the number of experimental observations, we cannot enumerate all probable networks at arbitrary moment, and calculate the average degree distribution directly. It has been proved that SPR is efficient in process networks with nodes addition, so a case study is provided here to confirm that the ESPR can process nodes preferential deletion. Assume a network consisted of 4 nodes as shown in the Fig.\ref{fig.2}, in which each node may have different degrees, and they will be deleted with probability $q_v(t),v=1,2,3,4$ from $t$ to $t+1$. \\
\indent We set the preferential deletion probability as $q_v=\frac{k(t)}{\sum_{w}k_w(t)},v=1,2,3,4$, so for nodes denoted as $1$ to $4$ have probabilities of $\frac{1}{8}$, $\frac{3}{8}$, $\frac{2}{8}$, $\frac{2}{8}$ respectively. The matrix of state transition probabilities $P$ is given by\\
\begin{equation}
	\left(
	\begin{array}{ccc}
		\frac{75}{192} & \frac{110}{192} & \frac{7}{192}\\
		\frac{3}{96} & \frac{74}{96} & \frac{19}{96}\\
		\frac{3}{64} & \frac{14}{64} & \frac{47}{64}\\
	\end{array}
	\right)
\end{equation}
Further, the average state distribution for new network with 3 nodes at $t+1$ can be expressed as
\begin{align*}
	P_{(3,0)}(t+1)=\frac{1}{4}\times\frac{75}{192}+\frac{1}{2}\times\frac{3}{96}+\frac{1}{4}\times\frac{3}{64}=\frac{1}{8}
\end{align*}
\begin{align*}
	P_{(3,1)}(t+1)=\frac{1}{4}\times\frac{110}{192}+\frac{1}{2}\times\frac{74}{96}+\frac{1}{4}\times\frac{14}{64}=\frac{7}{12}
\end{align*}
\begin{align*}
	P_{(3,2)}(t+1)=\frac{1}{4}\times\frac{7}{192}+\frac{1}{2}\times\frac{19}{96}+\frac{1}{4}\times\frac{47}{64}=\frac{7}{24}
\end{align*}

So the average degree distribution solved by extended SPR is $(\frac{1}{8},\frac{14}{24},\frac{7}{24},0)$.

\indent Basing on the intuitive evolving process in Fig.2, it is easy to determine that the state distributions for 4 networks are $(0,0,1,0)$, $(\frac{1}{3},\frac{2}{3},0,0)$,$(0,\frac{2}{3},\frac{1}{3},0)$,$(0,\frac{1}{3},\frac{2}{3},0)$ with the probability of $\frac{1}{8}$,$\frac{3}{8}$,$\frac{2}{8}$ and $\frac{2}{8}$ respectively at $t+1$. With the definition above, the average degree distribution is  
\begin{equation}
	\begin{aligned}
		P_K(t+1)&=\frac{1}{8}\times(0,0,1,0)+\frac{3}{8}\times(\frac{1}{3},\frac{2}{3},0,0)+\frac{2}{8}\times(0,\frac{2}{3},\frac{1}{3},0)+\frac{2}{8}\times(0,\frac{1}{3},\frac{2}{3},0)\\
		&=(\frac{1}{8},\frac{14}{24},\frac{7}{24},0)
	\end{aligned}
\end{equation}
Combining with the principle of SPR, proof completed.$\hfill\blacksquare$\\
\indent Moreover, it can be proved that the original SPR for nodes random deletion is a special case of ESPR we proposed above, which demonstrates that the ESPR is a complete framework to discuss the steady-state degree distribution for both nodes addtion and deletion. \\

\textbf{Theorem 2}. For a random birth and death network(RBDN) $\left\{G(t), t\geq0\right\}$, suppose the $SP_k$ represents the steady degree distribution of $
\left\{G(t)\right\}$ obtained from SPR, the $EP_k$ represents the steady degree distribution of $\left\{G(t)\right\}$ obtained from ESPR. Then, the formula $SP_k=EP_k$ holds.\\
\textbf{Proof}. For a RBDN, if it performed node deletion operation at $t$, a node $v$ in the network possessing $n$ nodes with degree $k$ would be removed with probability $\frac{1}{n}$. According to the approaches we previously introduced, it means remove $v$ with $q_v(t)=\frac{1}{n}$ from the perspective of nodes preferential deletion. So the state transition probabilitis of $v$ in ESPR are as follows. 
 $(1)$Node $v$ with state $(n,k)$ is not isolated \\
\begin{equation}
	\begin{aligned}
		P_{(n,k)(n-1,k)}(t)&=\frac{n-1-k}{n}  \quad\quad\quad\quad\quad
		P_{(n,k)(n-1,k-1)}(t)&=\frac{k}{n}
	\end{aligned}
\end{equation}
$(2)$Node $v$ with state $(n,k)$ is isolated\\
\begin{equation}
	\begin{aligned} 
		P_{(n,k)(n-1,k^{'})}(t)&=q_v(t)[\frac{P_{k^{'}}(t)(1-q_{(n,k^{'})}(t))}{\sum_kP_{(n,k)}(t)(1-q_{(n,k)}(t))}\frac{P_{(n,k^{'})(n-1,k^{'})}(t)}{P_{(n,k^{'})(n-1,k^{'})}(t)+P_{(n,k^{'})(n-1,k^{'}-1)}(t)}\\
		& \quad +\frac{P_{k^{'}+1}(t)(1-q_{(n,k^{'}+1))}(t)}{\sum_kP_k(t)(1-q_{(n,k)}(t))}\frac{P_{(n,k^{'}+1)(n-1,k^{'})}(t)}{P_{(n,k^{'}+1)(n-1,k^{'}+1)}(t)+P_{(n,k^{'}+1)(n-1,k^{'})}(t)}]\\
		&=\frac{nq_v(t)}{n-1}[P_{k^{'}}(t)P_{(n;k^{'})(n-1,k^{'})}(t)+P_{k^{'}+1}(t)P_{(n,k^{'}+1)(n-1,k^{'})}(t)]\\
		&=\frac{1}{n-1}[P_{k^{'}}(t)\frac{n-1-k^{'}}{n}+P_{k^{'}+1}(t)\frac{k^{'}+1}{n}]
	\end{aligned}
\end{equation}	\\
So referring to the approach in section \ref{sec.3} we get the final state transition equation in ESPR shown as below.\\
\begin{equation}
	\begin{aligned}
		P_{(n-1,k)}(t+1)&=\sum_{k^{'}}^{\infty}P_{(n,k)}(t) P_{(n,k^{'})(n-1,k)}(t)\\
		&=P_{(n,0)}(t) P_{(n,0)(n-1,k)}(t)+P_{(n,1)}(t) P_{(n,1)(n-1,k)}(t)+P_{(n,2)}(t) P_{(n,2)(n-1,k)}(t)+...\\
		&\quad +P_{(n,k)}(t) P_{(n,k+1)(n-1,k)}(t)+P_{(n,k+1)}(t) P_{(n,k+1)(n-1,k)}(t)+...+P_{(n,n)}(t) P_{(n,n)(n-1,n)}(t)\\
		&=[P_{(n,0)}(t)+P_{(n,1)}(t)+...+P_{(n,k-1)}(t)+P_{(n,k+2)}(t)+...+P_{(n,n)}(t)]\frac{1}{n-1}[P_{(n,k)}(t)\frac{n-k-1}{n}\\
		&\quad +P_{(n,k+1)}(t)\frac{k+1}{n}]+P_{(n,k)}(t){\frac{n-k-1}{n}+\frac{1}{n-1}[P_{(n,k)}(t)\frac{n-k-1}{n}+P_{(n,k+1)}(t)\frac{k+1}{n}]}\\
		&\quad +P_{(n,k+1)}(t){\frac{k+1}{n}+\frac{1}{n-1}[P_{(n,k)}(t)\frac{n-k-1}{n}+P_{(n,k+1)}(t)\frac{k+1}{n}]}\\
		&=[1-P_{(n,k)}(t)-P_{(n,k+1)}(t)]\frac{1}{n-1}[P_{(n,k)}(t)\frac{n-k-1}{n}+P_{(n,k+1)}(t)\frac{k+1}{n}]+P_{(n,k)}(t)\frac{n-k-1}{n}\\
		&\quad +P_{(n,k+1)}(t)\frac{k+1}{n}+P_{(n,k)}(t)\frac{1}{n-1}[P_{(n,k)}(t)\frac{n-k-1}{n}+P_{(n,k+1)}(t)\frac{k+1}{n}]\\
		&\quad +P_{(n,k+1)}(t)\frac{1}{n-1}[P_{(n,k)}(t)\frac{n-k-1}{n}+P_{(n.k+1)}(t)\frac{k+1}{n}]\\
		&=\frac{1}{n-1}[P_{(n,k)}(t)\frac{n-k-1}{n}+P_{(n,k+1)}(t)\frac{k+1}{n}]+P_{(n,k)}(t)\frac{n-k-1}{n}+P_{(n,k+1)}(t)\frac{k+1}{n}\\
		&=P_{(n,k)}(t)\frac{n-k-1}{n-1}+P_{(n,k+1)}(t)\frac{k+1}{n-1}\\
	\end{aligned}	
\end{equation}\\
Similarly, referring to Appendix \ref{sec.7}, we get the tate transition equation in SPR:\\ 
\begin{equation}
	\begin{aligned}
	P_{(n-1,k)}{(t+1)}&=P_{(n,k)}(t)P_{(n,k)(n-1,k)}+P_{(n,k+1)}(t)P_{(n,k+1)(n-1,k)}\\
	&=P_{(n,k)}(t)\frac{n-k-1}{n-1}+P_{(n,k+1)}(t)\frac{k+1}{n-1}	
	\end{aligned}	
\end{equation}

\indent It is obvious that the results calculated by ESPR and SPR are both identical, illustrating the validity of the theorem. Proof completed.$\hfill\blacksquare$\\

\section{Conclusion}\label{sec.5}

As one of the basic evolutionary mechanisms for evolving networks, nodes preferential deletion has the least discusssion comparing with other mechanisms. It can be attributed to the limitation of methods to handle related statistics of this kind of networks. In fact, the SPR is an approach which studies networks with exact solution of the steady degree distribution, keeping the topology and statistics of network during evolving process. For this method has obvious advantages in addressing those discrete problem, we discussed how to model the actual networks when node is removed with probability propotional to its degree with an extension of SPR in this paper. \\
\indent The extend SPR(ESPR) we proposed has been proved that it can handle  evolving networks with nodes preferential removal, through comparing with the results of the average degree distribution between evolving network in reality and networks in ESPR. Besides, the fact that original SPR is a special case of ESPR has also been proved, i.e. ESPR is a sophisticated frame in solving steady degree distribution for evolving networks. It should be noted that the ESPR may still be complex to nonhomogeneous evolving network process completely, such as nodes and edges attached on them are all deleted, but it can solve the networks that remove nodes but rewire these  associated edges. Moreover, this frame may be applied to discuss other statistics of networks e.g. degree-degree correlation further, for it possesses obvious advantage for discrete situation than continuum formalism. \\

\section*{Acknowledgment}
This work was supported by the National Natural Science Foundation of China (NO. 61273015).\\

\begin{appendices}
\section{Nodes Addition}\label{sec.6}
For nodes random addition in SPR, specific approaches are given as follows.\\

\indent\textbf{Step1}: Select $n+1$ networks with the same topology from $G_n(t),n=1,2,...,N(t)$. Choose a network from them randomly and delete all edges of this network, then put these $n$ networks left after former approach into $A(t)$ and $n$ isolated nodes into the set $B(t)$. Take the same operation to the rest networks in $G_n(t)$.\\
\indent\textbf{Step2}: Select a node and network from $B(t)$ and $A(t)$ respectively. Connect the node to the network with a given mechanism, random or preferential attachment, and put this new network with $n+1$ nodes into $G_{n+1}(t+1)$ . Finally, repeat the same operation to the rest members in $A(t)$ and $B(t)$.\\

The specific approaches of nodes addition in SPR are shown in Fig.3. Thus, for a node having a state of $(n,k)$ at $t$, may change to $(n+1,k)$ or $(n+1,k+1)$ at $t+1$ after growth of network, and transition probabilities for this node is shown as follows. \\
\begin{align}
	\begin{split}
		P_{(n,k)(n+1,k)}=\frac{(n-m)p}{n+1}\quad\quad\quad
		P_{(n,k)(n+1,k+1)}=\frac{mp}{n+1}\quad\quad\quad
		P_{(n,k)(n+1,m)}=\frac{p}{n+1}
	\end{split}
\end{align}
Here, $m$ represents the number of nodes to which the new added node attaches.\\

\begin{figure}[!t]
	\centering
	\includegraphics[width=5in]{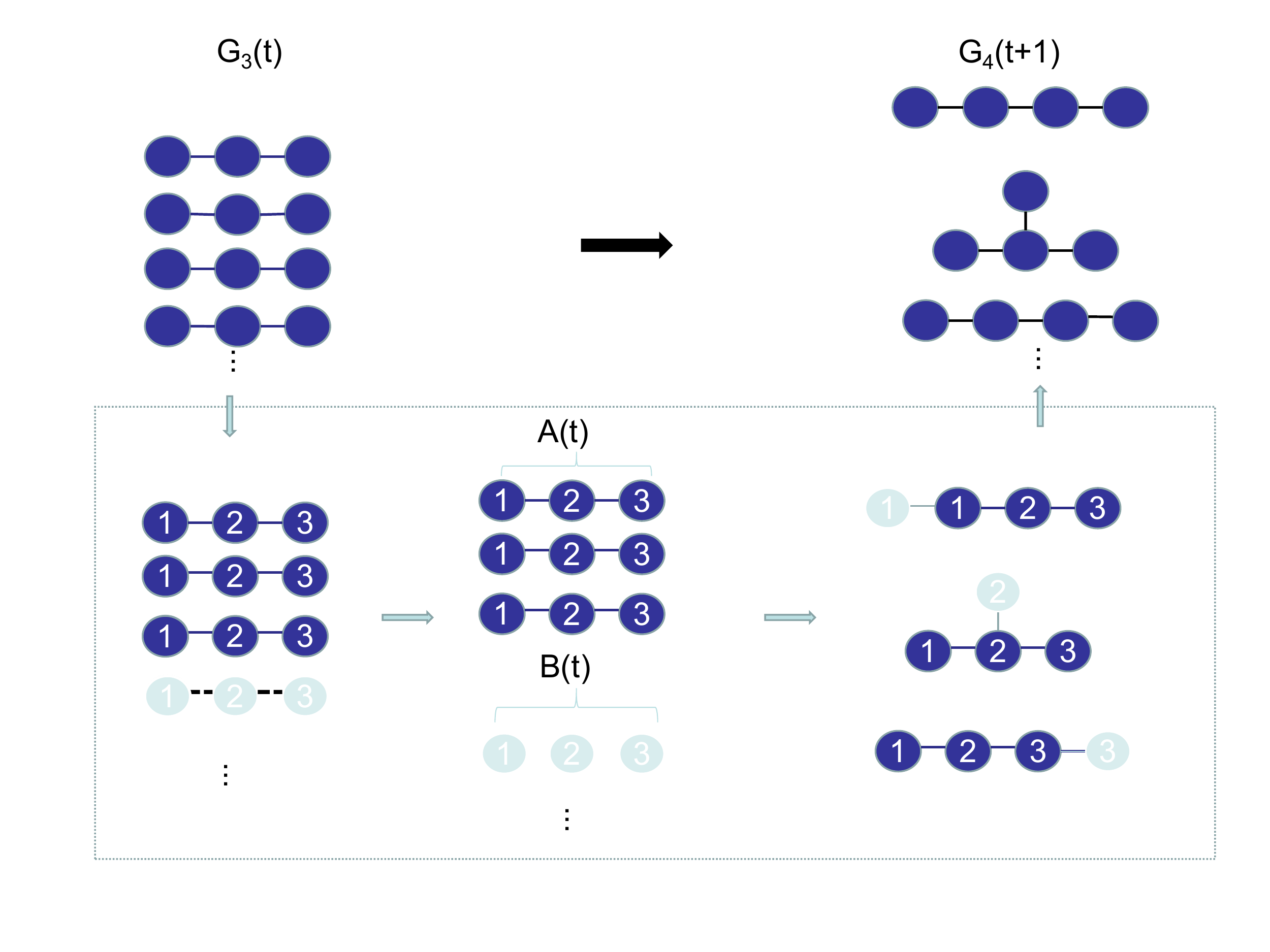}
	\caption{nodes random addition.}
	\label{fig.3}	
\end{figure}

\section{Nodes Random Deletion}\label{sec.7}
For nodes random deletion in SPR, it is more complicated than nodes addition introduced above and specific approaches are given as follows. \\

\indent\textbf{Step1}: $N$ networks with the same topology are selected from $G_n(t),n=1,2,...,i,...,N_t$, where $N$ is a common denominator of $n(n-1)$. Choose a network from them and delete a node from this network randomly. The network with $n-1$ nodes is placed into $A(t)$, the isolated node is placed into $B(t)$. Repeat the same operation to the rest networks in $G_n(t)$. After nodes removing, we get $B(t)=\cup_iB_i(t)$, where $B_i(t)$ represents the set of nodes at $t$ with different types in their original networks respectively, we distinguish them with notation of $i$.\\
\indent\textbf{Step2}: Select a network from $A(t)$ randomly, and select $n-1$ deleted nodes from $B(t)$ with the same nodes types as that of the chosen network. Connect these nodes with edges, referring to the topology of chosen network. Repeat the same operation to the rest members in $A(t)$ and $B(t)$. \\
\indent\textbf{Step3}: Put these new constructed networks into $C(t)$, and all networks contained in $A(t)$ and $C(t)$ are put into $G_{n-1}(t+1)$.\\

\indent In Fig.4, the details of operations in SPR are shown, especially the reconstruction of isolated nodes during the network reduction with specific numerical denotions for nodes as follows. 
\begin{figure}[!t]
	\centering
	\includegraphics[width=5in]{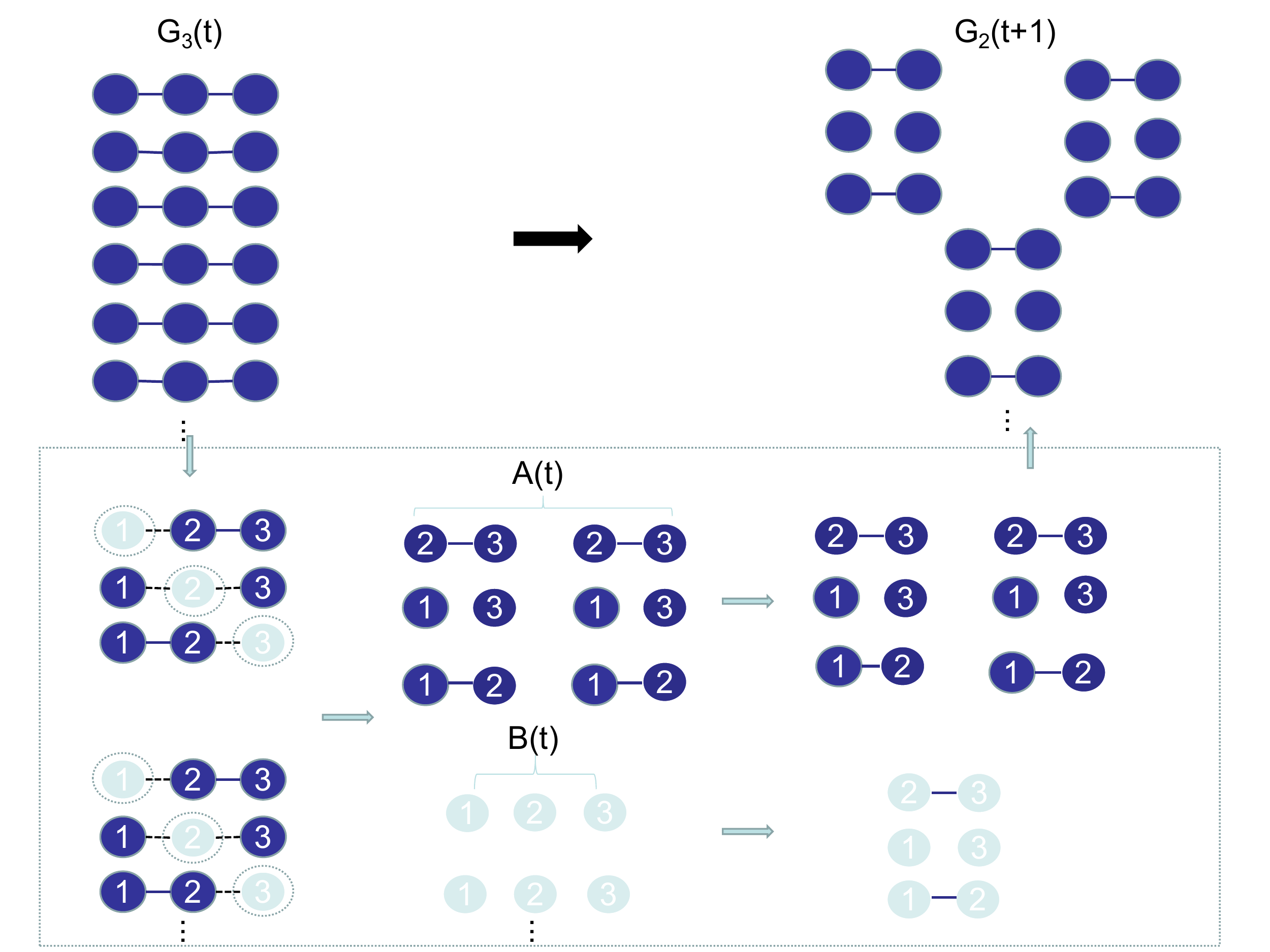}
	\caption{nodes random deletion.}
	\label{fig.4}
\end{figure}
A node with the state of $(n,k)$ at $t$ would get new state of $(n-1,k)$ or $(n-1,k-1)$ at $t+1$, when in the network reduction. So the transition probabilities are given as below.\\
\begin{align}
	\begin{split}
		P_{(n,k)(n-1,k)}=\frac{(n-k-1)q}{n-1} \quad\quad\quad\quad
		P_{(n,k)(n-1,k-1)}=\frac{kq}{n-1}
	\end{split}
\end{align}\\
\end{appendices}

\bibliography{paper}
\end{document}